\documentclass[12pt]{article}
\usepackage{amsmath, amsthm, amssymb, stmaryrd, fullpage}
\usepackage[all]{xy}

\newtheorem{theorem}{Theorem}[subsection]
\numberwithin{equation}{theorem}
\newtheorem{lemma}[theorem]{Lemma}
\newtheorem{conj}[theorem]{Conjecture}
\newtheorem{cor}[theorem]{Corollary}
\newtheorem{prop}[theorem]{Proposition}

\theoremstyle{definition}
\newtheorem{defn}[theorem]{Definition}

\newtheorem{remark}[theorem]{Remark}
\newtheorem{convention}[theorem]{Convention}

\def\AAA{\mathbb{A}}

\def\PP{\mathbb{P}}
\def\QQ{\mathbb{Q}}
\def\RR{\mathbb{R}}
\def\ZZ{\mathbb{Z}}
\newcommand{\calE}{\mathcal{E}}
\newcommand{\calF}{\mathcal{F}}
\newcommand{\calG}{\mathcal{G}}
\def\calO{\mathcal{O}}
\newcommand{\gothm}{\mathfrak{m}}
\newcommand{\gotho}{\mathfrak{o}}
\newcommand{\gothp}{\mathfrak{p}}

\def\alg{\mathrm{alg}}

\def\bv{\mathbf{v}}

\DeclareMathOperator{\Frac}{Frac}
\DeclareMathOperator{\height}{height}
\DeclareMathOperator{\Hom}{Hom}
\DeclareMathOperator{\Isoc}{Isoc}
\DeclareMathOperator{\rank}{rank}
\DeclareMathOperator{\ratrank}{ratrank}
\DeclareMathOperator{\Spec}{Spec}
\DeclareMathOperator{\trdeg}{trdeg}
\DeclareMathOperator{\unr}{unr}

\newcounter{fixmectr}

\begin{document}

\title{Semistable reduction for overconvergent $F$-isocrystals, II:
A valuation-theoretic approach}
\author{Kiran S. Kedlaya \\ Department of Mathematics \\ Massachusetts
Institute of Technology \\ 77 Massachusetts Avenue \\
Cambridge, MA 02139 \\
\texttt{kedlaya@mit.edu}}
\date{September 2, 2007}

\maketitle

\begin{abstract}
We introduce a valuation-theoretic approach to the problem of semistable
reduction (i.e., existence of logarithmic extensions on suitable covers)
of overconvergent isocrystals with Frobenius structure. The key tool
is the quasicompactness of the Riemann-Zariski space associated to the 
function field of a variety.
We also make some initial reductions, which allow attention to be focused
on valuations of height 1 and transcendence degree 0.
\end{abstract}

\tableofcontents

\section{Introduction}

This paper is the second of a series starting with \cite{kedlaya-part1}.
The goal of the series is to prove a 
``semistable reduction'' theorem for overconvergent $F$-isocrystals,
a class of $p$-adic analytic objects associated to schemes of finite type
over a field of characteristic $p>0$.
Such a theorem is expected to have consequences for the theory of rigid
cohomology, in which overconvergent $F$-isocrystals play the role of
coefficient objects.

In \cite{kedlaya-part1}, it was shown that the problem of extending
an overconvergent isocrystal on a variety $X$ to a log-isocrystal
on a larger variety $\overline{X}$ is governed by the triviality of some sort of
``local monodromy'' along components of the complement of $X$.
In this paper, we 
give a valuation-theoretic interpretation of this concept, which suggests an 
approach to the semistable reduction problem
to be pursued later in this series.

The context of this result (including a complex analogue) and a description
of potential applications is already given in the introduction of 
\cite{kedlaya-part1}, so we will not repeat it here. Instead, we devote
the remainder of this introduction to an overview of the results
specific to this paper, and a survey of the structure of the
various chapters of the paper.

\subsection{Valuations and semistable reduction}

Let $X \hookrightarrow \overline{X}$ be an open immersion of
varieties over a field
$k$ of characteristic $p>0$, with $X$ smooth and $\overline{X}$ proper,
and let $\calE$ be an $F$-isocrystal
(isocrystal with Frobenius structure)
on $X$ overconvergent
along $Z = \overline{X} \setminus X$. The semistable reduction problem,
as described in \cite[Section~7]{kedlaya-part1}, is to show that
$\calE$ admits a logarithmic extension with nilpotent residues after
being pulled back along some generically finite cover of $\overline{X}$.
When $X$ is a curve, this can be deduced from the $p$-adic local monodromy
theorem ($p$LMT) of Andr\'e \cite{andre}, Mebkhout \cite{mebkhout}, and the present
author \cite{kedlaya-local}. This derivation is carried out in \cite{kedlaya-semicurve};
the main point is that one can work locally, constructing the logarithmic
extension separately for each point of $Z$.

When $X$ has dimension greater than 1, one can still apply the $p$LMT
along codimension 1 components of $Z$, but one only obtains
a result that holds after ignoring a proper closed subset of the component. This
would be fine if one were always able to use a finite cover in the pullback
(by the analogues of Zariski-Nagata purity derived in \cite{kedlaya-part1}),
but that is not always possible: the result may be forced not to be smooth, in which
case some blowing up is required, producing additional components of codimension 1 along which it is not clear that any control on monodromy has been imposed.

To get around this, it is helpful to think of the application of the $p$LMT
as being parametrized by divisorial valuations, i.e., certain points
on the Riemann-Zariski space associated to the function field of $X$. One
is then naturally led to 
propose a version of the semistable reduction problem which
is local in Riemann-Zariski space. In this paper, we formulate the local
semistable reduction problem and explain its equivalence to the original
semistable reduction problem, 
using the quasicompactness of the Riemann-Zariski space. We also perform
some simplifying reductions that allow us to focus on what we call
\emph{minimal} valuations.
We defer a direct assault on the
local semistable reduction problem to subsequent papers in this series.

\subsection{Structure of the paper}

We conclude this introduction with a summary of the structure of the
paper.

In Section~\ref{sec:valuations}, we review some relevant 
facts from valuation theory, most notably the construction of Riemann-Zariski spaces.

In Section~\ref{sec:setup}, we describe the valuation-theoretic setup in 
more detail, formulating a local semistable reduction problem and verifying
that it is equivalent to the semistable reduction problem described
in \cite[Section~7]{kedlaya-part1}.

In Section~\ref{sec:reduction}, we show that the local semistable reduction
problem can be somewhat simplified. Specifically, we show that it suffices to
solve it when $k$ is algebraically closed, and 
the center valuation is of height 1 and has residue field $k$.

\subsection*{Acknowledgments}

Thanks to  Francesco Baldassarri and Pierre Berthelot for organizing
a useful workshop on $F$-isocrystals and rigid cohomology in 
December 2004/June 2005. 
Thanks to Daniel Caro for suggesting to consider partial compactifications
in detail.
The author was supported by NSF grant DMS-0400727.

\section{Review of valuation theory}
\label{sec:valuations}

In this chapter, we review some relevant facts from valuation theory,
notably the definition of the Riemann-Zariski space of a field.
We use the summary by Vaqui\'e \cite{vaquie} as our 
primary reference; in turn, Vaqui\'e's
 underlying primary reference is Zariski-Samuel \cite{zariski-samuel2}.

\setcounter{equation}{0}
\begin{convention}
For $A$ a local ring, let $\gothm_A$ denote the maximal ideal of $A$,
and let $\kappa_A = A/\gothm_A$ denote the residue field of $A$.
\end{convention}

\begin{convention}
By a \emph{variety} 
over a field $k$, we mean a reduced separated (but not necessarily irreducible) scheme of finite type over $k$. For $X$ an irreducible variety over $k$, let $k(X)$ 
denote the function field of $X$ over $k$.
By a \emph{smooth pair} over a field $k$, 
we mean a pair $(X,Z)$ consisting of a smooth 
irreducible $k$-variety $X$
and a strict normal crossings divisor $Z$ on $X$; we will conflate 
this pair with the log-scheme it determines.
\end{convention}

\subsection{Totally ordered groups}

We start with some standard facts about totally ordered groups,
which are used to discuss valuations.

\begin{defn}
By a \emph{totally ordered group}, we will mean an abelian group
$\Gamma$ written additively, equipped
with a total ordering $\leq$ with the property that
for $a,b,c \in \Gamma$,
$a \leq b$ if and only if $a+c \leq b+c$; as usual, we write
$a < b$ to mean that $a \leq b$ but $a \neq b$ (so that $b \not\leq a$),
and we write $a \geq b$ and $a > b$ to mean $b \leq a$ and $b < a$, respectively.
Write 
$\Gamma_\infty$ for the monoid $\Gamma \cup \{\infty\}$ in 
which $x + \infty = \infty$ for all $x \in \Gamma_\infty$, and extend
the total ordering to $\Gamma_\infty$ by declaring that for all $x \in \Gamma$,
$x < \infty$.
\end{defn}

\begin{defn}
Let $\Gamma$ be a totally ordered group. 
A proper subgroup $\Delta$ of $\Gamma$ is called an \emph{isolated subgroup}
if for any $\alpha \in \Delta$, $\beta \in \Gamma$ with
$\alpha \geq \beta \geq 0$, it follows that $\beta \in \Delta$.
It is easily shown that the isolated subgroups are totally ordered
under inclusion; define the
\emph{rank} of $\Gamma$, denoted $\rank(\Gamma)$, to be the
cardinality of the set of isolated 
subgroups of $\Gamma$.
\end{defn}

\begin{prop} \label{P:archimedean}
Let $\Gamma$ be a totally ordered group. Then the following are equivalent.
\begin{enumerate}
\item[(a)]
$\Gamma$ is isomorphic, as a totally ordered group,
to a subgroup of $\RR$ with its usual ordering.
\item[(b)]
$\rank(\Gamma) \leq 1$.
\item[(c)]
$\Gamma$ is archimedean: that is, for any $a,b \in \Gamma$ with
$b > 0$, there exists an integer $n$ with
$nb \geq a$.
\end{enumerate}
\end{prop}
\begin{proof}
See \cite[\S VI.10, p.\ 45]{zariski-samuel2}.
\end{proof}
\begin{cor}
A totally ordered group $\Gamma$ has rank at most $n$ if and only if it
is isomorphic to a subgroup of $\RR^n$ with its lexicographic ordering.
\end{cor}

\begin{defn}
Let $\Gamma$ be a totally ordered group.
Define the 
\emph{rational rank} of $\Gamma$, denoted $\ratrank(\Gamma)$, to be the
dimension of the $\QQ$-vector space $\Gamma \otimes_{\ZZ} \QQ$.
\end{defn}

\begin{prop} \label{P:rank rat}
Let $\Gamma$ be a totally ordered group, let
$\overline{\Gamma}$ be an isolated subgroup of $\Gamma$, 
and equip $\Gamma' = \Gamma/\overline{\Gamma}$ with the induced total ordering.
Then
\begin{align*}
\ratrank(\Gamma) &= \ratrank(\overline{\Gamma}) + \ratrank(\Gamma') \\
\rank(\Gamma) &= \rank(\overline{\Gamma}) + \rank(\Gamma') \\
\rank(\Gamma) &\leq \rank(\overline{\Gamma}) + \ratrank(\Gamma').
\end{align*}
In particular, $\rank(\Gamma) \leq \ratrank(\Gamma)$.
\end{prop}
\begin{proof}
See \cite[Proposition~3.5]{vaquie}.
\end{proof}

\subsection{Valuations}

\begin{defn} \label{D:val}
Let $R$ be a ring and let $\Gamma$ be a totally ordered group.
A \emph{valuation} (or \emph{Krull valuation})
on $R$ with values in $\Gamma$ 
is a function $v: R \to \Gamma_\infty$ satisfying the
following conditions.
\begin{enumerate}
\item[(a)] For $x,y \in R$, $v(xy) = v(x) + v(y)$.
\item[(b)] For $x,y \in R$, $v(x+y) \geq \min\{v(x), v(y)\}$.
\item[(c)] We have $v(1) = 0$ and $v(0) = \infty$.
\end{enumerate}
We write \emph{real valuation} as shorthand for
``Krull valuation with values in $\RR$''.
We say that the valuations $v_1, v_2$ on $R$, with
values in $\Gamma_1, \Gamma_2$, are \emph{equivalent} if
for all $x,y \in R$,
\[
v_1(x) \geq v_1(y) \qquad \Longleftrightarrow \qquad
v_2(x) \geq v_2(y).
\]
Define the \emph{value group} $\Gamma_v$ of $v$ to be the image of
$v$, as a totally ordered group; equivalent valuations have
isomorphic value groups, and in fact every valuation is equivalent to
a surjective valuation (in which $\Gamma = \Gamma_v$).
\end{defn}

\begin{remark}
With notation as in Definition~\ref{D:val},
the inverse image $\gothp = v^{-1}(\infty)$ is a prime ideal of $R$,
and the induced map $v: R/\gothp \to \Gamma_\infty$ extends uniquely
to a valuation on $\Frac(R/\gothp)$ \cite[Proposition~2.2]{vaquie}.
\end{remark}

\begin{defn}
For $F$ a field and $v: F \to \Gamma_\infty$ a valuation on $F$,
the subset $\gotho_v = \{x \in F: v(x) \geq 0\}$ is a valuation ring of $F$
(a subring of $F$ maximal among local rings with fraction field $F$
under domination; see \cite[Th\'eor\`eme~1.1]{vaquie} for other
characterizations). In this case,
we write $\kappa_v$ for the residue field $\kappa_{\gotho_v}$.
Conversely, every valuation ring of $F$ arises
from a unique equivalence class of valuations
\cite[Proposition~2.3]{vaquie}. 
Let $\gothm_v = \{x \in F: v(x) > 0\}$ denote the maximal ideal of
$\gotho_v$.
\end{defn}

\begin{defn}
For $k$ an integral domain, there is a valuation $v: k \to \{0, \infty\}$
with $v(x) = \infty$ if and only if $x=0$. We call $v$ the
\emph{trivial valuation} (or \emph{improper valuation}) on $k$.
More generally, if $R$ is a $k$-algebra, a \emph{valuation on 
$R$ over $k$} is a valuation which restricts to the trivial
valuation on $k$.
\end{defn}

\begin{defn}
Define the \emph{height} (or \emph{rank}, or \emph{real rank})
and \emph{rational rank} of $v$ as the rank and rational rank,
respectively, of the value group of $v$, and denote these by
$\height(v)$ (or $\rank(v)$) and $\ratrank(v)$, respectively.
By Proposition~\ref{P:rank rat}, we have
$\height(v) \leq \ratrank(v)$.
Also, note that by Proposition~\ref{P:archimedean}, a valuation has height 1
if and only if it is equivalent to a real valuation.
\end{defn}

Valuations of height greater than 1 can be written as ``compositions''
of valuations of smaller height.
\begin{defn} \label{D:composition}
Let $F$ be a field, and let $v: F \to \Gamma_\infty$ be a valuation.
For $\overline{\Gamma}$ a nontrivial isolated subgroup of $\Gamma_v = 
v(F)$, put
$\Gamma' = \Gamma/\overline{\Gamma}$, and let
$v': F \to \Gamma'_\infty$ be the composition of $v$ with
the quotient map $\Gamma \to \Gamma'$; it is again a valuation.
Let $\overline{v}: \kappa_{v'} \to \overline{\Gamma}_\infty$ be the
map induced by $v$; it too is a valuation. In this situation, we
write $v = v' \circ \overline{v}$ and say that
$v$ is the \emph{composition} of $v'$ and $\overline{v}$; note 
that \cite[Proposition~3.5]{vaquie}
\begin{align*}
\height(v) &= \height(v') + \height(\overline{v}) \\
\ratrank(v) &= \ratrank(v') + \ratrank(\overline{v}).
\end{align*}
\end{defn}

\begin{remark}
Note that the convention ``overbars denote reduction'' here
applies to the valuation rings, rather than to the value groups.
\end{remark}

\begin{defn} \label{D:extension}
Let $E/F$ be an extension of fields. If $w$ is a valuation on $E$, the
restriction of $w$ to $F$ is a valuation on $F$; if that valuation is $v$,
we say that $w$ is an \emph{extension to $E$} of $v$. Note that this happens if
and only if $\gotho_w \cap F = \gotho_v$; in particular, since every
integral domain is contained in a valuation ring (by a Zornication), every
valuation on $F$ admits at least one extension to $E$.
\end{defn}

\begin{defn}
Let $E/F$ be a finite extension of fields, let $v$ be a valuation on $F$, and let $w$ be an extension to $E$ of $v$. Define the \emph{ramification index} 
\[
e(w/v) = [\Gamma_w:\Gamma_v]
\]
and the \emph{residual degree}
\[
f(w/v) = [\kappa_w:\kappa_v].
\]
These numbers are both finite and satisfy the fundamental inequality
\[
e(w/v) f(w/v) \leq [E:F]
\]
\cite[Proposition~5.1]{vaquie}. In particular,
$\rank(w) = \rank(v)$ and $\ratrank(w) = \ratrank(v)$.
\end{defn}

\subsection{Riemann-Zariski spaces}

We now recall the definition of a Riemann-Zariski space,
following \cite[\S 7]{vaquie}.

\begin{defn}
Let $F$ be a field and let $R$ be a subring of $F$. Define the 
\emph{Riemann-Zariski space} $S_{F/R}$ 
as the set consisting of the equivalence classes of
valuations on $F$ which are nonnegative on $R$.
This set carries two natural topologies, the coarser
\emph{Zariski topology}
and the finer
\emph{patch topology} (or \emph{Zariski-Hausdorff topology}), in which bases
are given by sets of one of the following respective forms:
\begin{align*}
\mbox{Zariski}: \ &
\{
v \in S_{F/R}: v(x_1) \geq 0, \dots, v(x_m) \geq 0
\} \\
\mbox{patch}: \ &
\{
v \in S_{F/R}: v(x_1) \geq 0, \dots, v(x_m) \geq 0; \quad v(y_1) > 0,
\dots, v(y_n) > 0
\}
\end{align*}
for $x_1, \dots, x_m, y_1, \dots, y_n \in F$.
The patch topology is Hausdorff, while the Zariski
topology is only $T_1$ in general.
Note that if $E/F$ is a field extension, then there is a natural surjection
$S_{E/R} \to S_{F/R}$ obtained by restricting valuations from $E$ to $F$
(see Definition~\ref{D:extension} for the surjectivity),
which is continuous for either consistent choice of topologies.
\end{defn}

\begin{defn}
For $F$ a field, let $T_F$ denote the set of functions from $F$ to $\{0,+,-\}$,
equipped with the product topology associated to the discrete topology on
$\{0,+,-\}$. Consider the map from $S_{F/R}$ to $T_F$ that associates to
$v \in S_{F/R}$ the function $f_v \in T_F$ taking $x \in F$ to $0,+,-$ depending on whether
$v(x) = 0$, $v(x) > 0$, $v(x) < 0$. This map $S_{F/R} \to T_F$
is injective because from $f_v$
we can recover $\gotho_v$ as $\{x \in F: f_v \in \{0,+\}\}$, and hence can recover $v$.
Moreover, the subspace topology induced on $S_{F/R}$ is visibly the same as the
patch topology. 
We may similarly recover the Zariski topology by starting
with the topology on $\{0,+,-\}$ with open sets $\emptyset, \{0,+\},
\{0,+,-\}$.
\end{defn}

The fundamental property of Riemann-Zariski spaces
is the following result (see the proofs of \cite[Th\'eor\`eme~7.2]{vaquie},
\cite[Theorem~VI.40]{zariski-samuel2}).
\begin{theorem} \label{T:Riemann-Zariski}
Let $F$ be a field and let $R$ be a subring of $F$. 
Then $S_{F/R}$ with the patch topology is a closed subset of 
$T_F$. Consequently (by Tykhonov's theorem) $S_{F/R}$ is compact under the
patch topology and quasicompact under the Zariski topology.
\end{theorem}

\subsection{Centers of valuations}

\begin{defn}
Let $F$ be a field, let $v: F \to \Gamma_\infty$ be a valuation on
$F$, and let $A$ be a subring of $\gotho_v$. The \emph{center} of
$v$ on $A$ is the ideal $A \cap \gothm_v$ in $A$; it is prime because
it is the contraction of a prime ideal of $\gotho_v$. If $A$ is a local ring
and $F = \Frac(A)$,
we say $v$ is \emph{centered in $A$} if the center of $v$ on $A$ equals
$\gothm_A$; it is equivalent to say that the valuation ring of $v$
dominates $A$.
\end{defn}

The concept of the center of a valuation also extends to schemes.
\begin{defn} \label{D:center scheme}
Let $X$ be an integral scheme, and let $v$ be a valuation on the function
field of $X$. Then
the set of points $x$ whose local rings
$\calO_{X,x}$ are contained in the valuation ring $\gotho_v$ is an
irreducible or empty
closed subset of $X$ \cite[Proposition~6.2]{vaquie}; we call this set 
(or the corresponding reduced closed subscheme)
the \emph{center} of $v$ on $X$. Note that the generic point of the
center of $v$ is the unique point $x$ such that $v$ is centered in 
$\calO_{X,x}$.
If the center of $v$ on $X$ is nonempty, we say that $v$ is 
\emph{centered on $X$}. If $X$ is proper over a field, then $v$ is
always centered on $X$ \cite[Proposition~6.3]{vaquie}. (See also 
Lemma~\ref{L:rel center} below.)
\end{defn}

\begin{prop} \label{P:center dim}
Let $X$ be a proper irreducible variety over a field $k$,
and let $v$ be a valuation on $k(X)$ over $k$. 
Then the dimension of the center of $v$ on $X$ is at most 
$\trdeg(\kappa_v/k)$. Moreover, there exists a blowup $Y$ of $X$ such that 
the center of $v$ on $Y$ has dimension equal to
$\trdeg(\kappa_v/k)$, and the same is true of any further blowup $Y'$ of $Y$.
\end{prop}
\begin{proof}
The equality occurs whenever the local ring on $Y$ of the generic point of
the center of $v$ contains elements lifting a transcendence basis of 
$\kappa_v$ over $k$; this clearly occurs for some $Y$ and continues to occur after further blowing up. Compare \cite[Proposition~6.4]{vaquie}.
\end{proof}

Note the following relationship to the Riemann-Zariski space 
\cite[Proposition~7.4]{vaquie}.
\begin{prop} \label{P:open subset}
Let $X$ be an integral noetherian scheme over a field $k$.
Then the set of valuations $v \in S_{k(X)/k}$ with nonempty center on $X$ is
an open subset $U(X)$ of $S_{k(X)/k}$ for the Zariski topology, 
and the map $U(X) \to X$ carrying
$v \in U(X)$ to the generic point of its center is continuous for the Zariski
topology on $S_{k(X)/k}$ and the usual (Zariski) topology on $X$.
\end{prop}

\begin{lemma} \label{L:rel center}
Let $f: X' \to X$ be a proper morphism between irreducible varieties over $k$.
Let $w$ be a valuation on $k(X')$ over $k$, and let $v$ be the restriction of
$w$ to $k(X)$. Let $Y$ be the center of $v$ on $X$, and let $Y'$ be the center
of $w$ on $X'$. Then the generic point of $Y'$ maps to the generic
point of $Y$.
\end{lemma}
\begin{proof}
The proof is a relative version of \cite[Proposition~6.3]{vaquie}.
In the diagram
\[
\xymatrix{
\Spec k(X') \ar[r] \ar[d] & X' \ar^f[d] \\
\Spec \gotho_w \ar[r] \ar@{-->}[ru] & X,
}
\]
the lower horizontal arrow factors as the surjection
$\Spec \gotho_w \to \Spec \gotho_v$ followed by the map
$\Spec \gotho_v \to X$ with image $Y$.
By the valuative criterion for
properness, into the diagram we can insert a map $\Spec \gotho_w \to X'$,
whose image of the closed point of $\Spec \gotho_w$ is the generic point
of $Y'$. This proves the claim.
\end{proof}

\subsection{Abhyankar's inequality}

We recall a fundamental result of Abhyankar \cite[Th\'eor\`eme~9.2]{vaquie}.

\begin{defn}
We say a valuation $v$ is \emph{discrete} if
its value group is isomorphic to $\ZZ^d$ under the lexicographic
ordering, for some nonnegative integer $d$. Note that this is
more inclusive than the layman's definition of a ``discrete valuation'';
that concept corresponds in standard valuation-theoretic terminology
to a \emph{divisorial valuation}. See Definition~\ref{D:divmon}.
\end{defn}

\begin{theorem}[Abhyankar] \label{T:abhyankar}
Let $A$ be a noetherian local ring, and put $F = \Frac(A)$. Let
$v: F \to \Gamma_\infty$ be a valuation on $F$ centered in $A$.
\begin{enumerate}
\item[(a)]
The following inequality holds:
\begin{equation} \label{eq:abhyankar}
\ratrank(v) + \trdeg(\kappa_v/\kappa_A) \leq \dim(A).
\end{equation}
\item[(b)]
If equality holds in (a), then $\Gamma \cong \ZZ^{\ratrank(v)}$ (as an abstract
group) and $\kappa_v$ is a finitely generated field extension of $\kappa_A$.
\item[(c)]
If $\rank(v) + \trdeg(\kappa_v/\kappa_A) = \dim(A)$, then $v$ is discrete.
\end{enumerate}
\end{theorem}

\begin{defn} \label{D:divmon}
Let $A$ be a noetherian local ring, and put $F = \Frac(A)$.
Let $v: F \to \Gamma_\infty$ be a valuation on $F$ whose valuation ring contains $A$.
We say $v$ is an \emph{Abhyankar valuation} if equality holds in
\eqref{eq:abhyankar}. We say a valuation $v$ is \emph{divisorial} if
\[
\rank(v) = 1, \qquad \trdeg(\kappa_v/\kappa_A) = \dim(A)-1.
\]
We say that $v$ is
\emph{monomial} if
\[
\rank(v) = 1, \qquad \ratrank(v) = \dim(A), \qquad \kappa_v = \kappa_A;
\]
beware that some authors may prefer not to include the rank restriction.
\end{defn}

\begin{remark}
For $k$ a field, any valuation $v$ over $k$ 
on a finitely generated field $K$ over $k$ is subject to Abhyankar's
inequality. That is because we may choose a proper variety $X$ over
$k$ with function
field $K$, on which $v$ will be centered (see Definition~\ref{D:center
scheme}). In particular, $v$ is centered on some affine chart of $X$,
whose coordinate ring is noetherian, as is its localization at the 
center of $v$.
\end{remark}

\section{The local approach to semistable reduction}
\label{sec:setup}

In this chapter, we recall the statement of the 
semistable reduction problem, 
then demonstrate its equivalence with a formally weaker form 
which is in some sense local at a valuation. 

\setcounter{theorem}{0}
\begin{convention}
Throughout this chapter, let $k$ be a field of characteristic $p>0$,
and fix a power $q$ of the prime $p$.
Let $K$ be a discretely valued field of characteristic $0$ with residue field
$k$. Assume that there exists a continuous endomorphism $\sigma_K: K \to K$
lifting the $q$-power Frobenius morphism on $k$, and fix a choice of $\sigma_K$.
\end{convention}

\begin{remark}
The restriction to $K$ discretely valued is necessitated in part by
that restriction in the work of Shiho \cite{shiho1, shiho2} invoked in \cite{kedlaya-part1}, and in part by
that restriction in the construction of slope filtrations 
for Frobenius modules \cite{kedlaya-local}, which will intervene at more than
one point in this series.
\end{remark}

We retain notation and terminology as set in \cite{kedlaya-part1};
for convenience, we recall a few of the less standard notations.

\begin{defn}
Let $\Gamma^*$ be the divisible closure of $|K^*|$. We say an interval
$I \subseteq [0,+\infty)$ 
is \emph{aligned} if each endpoint at which $I$ is closed
is either zero or an element of $\Gamma^*$. For $I$ an aligned interval,
let $A^m_K(I)$ denote the admissible subspace of the rigid analytic $m$-plane
with coordinates $t_1, \dots, t_m$, consisting of points where
$|t_i| \in I$ for $i=1, \dots,m$. We drop the parentheses around $I$ when
it is written out explicitly; for instance, we write $A^m_K[0,1)$ for the open
unit polydisc. If we need to specify the name $t$ of the family of dummy 
variables, we add it to the notation as a subscript, e.g, $A^m_{K,t}(I)$.
\end{defn}

\subsection{Alterations}

Before proceeding further, we recall the fundamental concept of
alterations, from \cite[2.20]{dejong}.

\begin{defn}
For $X$ a noetherian integral scheme, an \emph{alteration} of $X$
is a proper dominant morphism $f: X_1 \to X$ with $X_1$ irreducible
and $f$ generically finite,
i.e., there exists $U \subseteq X$ open dense such that 
$f^{-1}(U) \to U$ is finite. 
\end{defn}

\begin{defn}
Let $X$ be an irreducible $k$-variety,
and let $Z$ be a proper closed subset of $X$. A \emph{quasiresolution} of the
pair $(X,Z)$ consists of
\begin{itemize}
\item
an alteration $f: X_1 \to X$ over $k$, which is required to be 
generically \'etale if $k$ is perfect, and
\item
an open immersion $j: X_1 \hookrightarrow
\overline{X_1}$ over $k$, with $\overline{X_1}$ projective over $k$,
\end{itemize}
such that
$(\overline{X_1}, j(f^{-1}(Z)) \cup (\overline{X_1} \setminus j(X_1)))$ 
form a smooth pair.
By a quasiresolution of $X$, we mean a quasiresolution of the pair $(X, \emptyset)$.
\end{defn}

In terms of this definition, de Jong's alterations theorem is as follows
\cite[Theorem~4.1]{dejong}.
\begin{theorem} \label{T:alterations}
Let $X$ be an irreducible $k$-variety,
and let $Z$ be a proper closed subset of $X$. Then the pair $(X,Z)$ admits a quasiresolution.
\end{theorem}

\subsection{The semistable reduction problem}

We now formulate the semistable reduction problem, following 
Shiho \cite[Conjecture~3.1.8]{shiho2}.
(This was done once already in \cite[Conjecture~7.1.2]{kedlaya-part1}.)

\begin{defn}
Let $(X,Z)$ be a smooth pair, and let $\calE$ be a convergent isocrystal on $U = X
\setminus Z$. We say that $\calE$ is \emph{log-extendable
(on $X$)} if $\calE$ extends to a convergent 
log-isocrystal
$\calF$ with nilpotent residues on the log-scheme $(X,Z)$, in the sense of
Shiho \cite{shiho1, shiho2}. 
\end{defn}

The following alternate formulation is 
\cite[Theorem~6.4.5]{kedlaya-part1} (the main result of \cite{kedlaya-part1}).
\begin{theorem} \label{T:log extend}
Let $(X,Z)$ be a smooth pair, and let $\calE$ be an isocrystal on
$U = X \setminus Z$ overconvergent along $Z$. Then $\calE$ has unipotent
monodromy along $Z$ in the sense of \cite[Definition~4.4.2]{kedlaya-part1}
if and only if $\calE$ is log-extendable on $X$. Moreover, the restriction
functor, from convergent log-isocrystals with nilpotent residues on $(X,Z)$
to isocrystals on $U$ overconvergent along $Z$, is fully faithful; in particular,
a log-extension with nilpotent residues is unique if it exists,
and any Frobenius on $\calE$ also acts on such an extension.
\end{theorem}

\begin{remark} \label{R:pullback log}
Let $(X,Z)$ and $(X',Z')$ be smooth pairs, and let
$f: X' \to X$ be a morphism with $f(X' \setminus Z') \subseteq X \setminus Z$.
If $\calE$ is a convergent isocrystal on $X \setminus Z$
which is log-extendable on $X$, then 
$f^* \calE$ is log-extendable
on $X'$; this follows from the existence of pullbacks in Shiho's category of 
convergent log-isocrystals.
\end{remark}

\begin{defn} \label{D:semi}
Let $X$ be a smooth irreducible $k$-variety, and
let $\calE$ be an overconvergent $F$-isocrystal on $X$. We say that 
$\calE$ \emph{admits semistable reduction} if there
exists a quasiresolution $(f: X_1 \to X, j: X_1 \hookrightarrow 
\overline{X_1})$ of
$X$ such that $f^* \calE$ is log-extendable on $\overline{X_1}$.
\end{defn}

We now recall, in the present terminology, the formulation of Shiho's conjecture
\cite[Conjecture~3.1.8]{shiho2} given earlier in this series
\cite[Conjecture~7.1.2]{kedlaya-part1}.
\begin{conj} \label{conj:Shiho}
Let $X$ be a smooth irreducible 
$k$-variety and let $\calE$ be an overconvergent
$F$-isocrystal on $X$. Then $\calE$ admits semistable reduction.
\end{conj}
Actually, Shiho's original conjecture only required $k$ perfect; however,
the distinction between this and the general case is illusory.
In fact, one may even reduce to considering algebraically
closed base fields, as follows.
\begin{prop} \label{P:only perfect}
Let $X$ be a smooth irreducible 
$k$-variety, and let $\calE$ be an overconvergent
$F$-isocrystal on $X$. 
Let $K^{\unr}$ be the maximal unramified extension of $K$,
and let $K'$ be the completion of the direct limit $K^{\unr} 
\stackrel{\sigma_K}{\to} K^{\unr}
\stackrel{\sigma_K}{\to} \cdots$.
Let $\calE'$ be the overconvergent $F$-isocrystal on
$X' = X \times_k k^{\alg}$, with coefficients in $K'$, obtained by
base extension from $\calE$. If $\calE'$ admits semistable reduction,
then so does $\calE$.
\end{prop}
\begin{proof}
If $\calE'$ admits semistable reduction, then there exists a quasiresolution
$(f': X'_1 \to X', j: X'_1 \hookrightarrow \overline{X'_1})$
such that $(f')^* \calE'$ is log-extendable on $\overline{X'_1}$.
Since specifying the data of this quasiresolution only involves a finite 
number of elements of $k^{\alg}$, we can realize it over some finite
extension $k'$ of $k$. This means we can produce an alteration
$f_1: X_1 \to X$ such that $X_1 \times_k k'$ is a disjoint union of copies of
$X'_1$. Unfortunately, $X_1$ need not be smooth over $k$; however, 
if we construct a quasiresolution
$(f_2: X_2 \to X_1, j_2: X_2 \hookrightarrow \overline{X_2})$, 
then the base extension of $(f_2 \circ f_1)^* \calE'$ to $k^{\alg}$
is log-extendable on $\overline{X_2} \times_k k^{\alg}$. Since
local unipotence can be checked after a field extension 
\cite[Remark~3.4.4]{kedlaya-part1}, we may apply
Theorem~\ref{T:log extend} to deduce that $(f_2 \circ f_1)^* \calE$ is also
log-extendable on $\overline{X_2}$. Hence $\calE$ admits semistable reduction.
\end{proof}

\subsection{Local semistable reduction}

We next formulate a local version of semistable reduction, then relate it
to global
 semistable reduction via the quasicompactness of Riemann-Zariski spaces.
\begin{defn} \label{D:local semi}
Let $X$ be a smooth irreducible $k$-variety, and
let $\calE$ be an overconvergent $F$-isocrystal on $X$. 
For $v$ a valuation on the function field $k(X)$ over $k$, we say that
$\calE$ \emph{admits local semistable reduction} at $v$ if there 
exists a quasiresolution $(f: X_1 \to X, j: X_1 \hookrightarrow 
\overline{X_1})$ of $X$
such that $f^* \calE$ is log-extendable on some open subscheme of $\overline{X_1}$
on which each extension of $v$ to $k(X_1)$ is centered.
\end{defn}
\begin{remark}
If $g:Y \to X$ is any alteration, and $g^* \calE$ admits local
semistable reduction at every extension of $v$ to $k(Y)$, then $\calE$ admits
local semistable reduction at $v$.
\end{remark}

\begin{lemma} \label{L:open}
Let $X$ be a smooth irreducible $k$-variety,
let $(f: X_1 \to X, j: X_1 \hookrightarrow \overline{X_1})$ be a quasiresolution
of $X$, and let $U$ be an open subscheme of $\overline{X_1}$. 
Then the set of $v \in S_{k(X)/k}$, 
all of whose extensions to $k(X_1)$ are centered on $U$,
is an open subset of $S_{k(X)/k}$ for the patch topology.
\end{lemma}
\begin{proof}
Let $A$ be the set of valuations on $k(X_1)$ centered on $U$; then $A$ is
open in $S_{k(X_1)/k}$. Put $B = S_{k(X_1)/k} \setminus A$, which is thus closed
for the patch topology; since $S_{k(X_1)/k}$ is compact by
Theorem~\ref{T:Riemann-Zariski}, so then is $B$. Let $C$ be the image of $B$
under the restriction map $S_{k(X_1)/k} \to S_{k(X)/k}$; then $C$ is quasicompact
since it is the image of a quasicompact topological space under a continuous map.
Since $S_{k(X)/k}$ is Hausdorff under the patch topology, 
$C$ is Hausdorff, hence compact, hence closed. The set we are looking
for is the complement of $C$, so we are done.
\end{proof}

\begin{prop} \label{P:local to global}
Let $X$ be a smooth irreducible $k$-variety, and
let $\calE$ be an overconvergent $F$-isocrystal on $X$. 
Suppose that $\calE$ admits local semistable reduction at every valuation 
on $k(X)$ over $k$.
Then $\calE$ admits semistable reduction.
\end{prop}
\begin{proof}
Consider the Riemann-Zariski space $S_{k(X)/k}$ equipped with the 
patch topology. By hypothesis, for each
$v_i \in S_{k(X)/k}$, we may choose a quasiresolution $(f_i: X_i \to X, j_i: X_i \hookrightarrow
\overline{X_i})$ of $X$ such that $f_i^* \calE$ is log-extendable on some
open subscheme $U_i$ of $\overline{X_i}$
containing the center of each extension of $v_i$ to $k(X_i)$.
Let $B_i$ be the set of valuations $w \in S_{k(X)/k}$ each of whose extensions
to $k(X_i)$ is centered in $U_i$; by Lemma~\ref{L:open}, $B_i$ is an open neighborhood
of $v_i$ in $S_{k(X)/k}$.

By Theorem~\ref{T:Riemann-Zariski}, $S_{k(X)/k}$ is compact, so there
exist finitely many valuations $v_1, \dots, v_n \in S_{k(X)/k}$ such that
$B_{1} \cup \cdots \cup B_{n} = S_{k(X)/k}$. 
Apply Theorem~\ref{T:alterations} (to the closure 
in $\overline{X_1} \times_k \cdots \times_k \overline{X_n}$ 
of an irreducible component of $X_1 \times_X \cdots \times_X X_n$)
to produce a smooth pair
$(\overline{Y},E)$ with $\overline{Y}$ projective, admitting maps
$g_i: (\overline{Y},E) \to (\overline{X_i},\overline{X_i} \setminus
X_i)$ for $i=1, \dots, n$.
Let $Y$ be the inverse image of 
$X_1 \times_X \cdots \times_X X_n$ in $\overline{Y}$, so that
$Y$ is open dense in $\overline{Y}$ and the $g_i$ induce a projective map
$g: Y \to X$.
Then by Remark~\ref{R:pullback log},
$g^* \calE = (g_i \circ f_i)^* \calE$ is log-extendable on $g_i^{-1}(U_i)$
for each $i$.

Let $H$ be a component of $E$; then $H$ corresponds to a divisorial
valuation on $k(Y)$, whose restriction to $k(X)$ 
must lie in one of the $B_{i}$.
For any such $i$, $g_i^{-1}(U_i)$ meets $H$, so by the easy direction of
Theorem~\ref{T:log extend} plus
\cite[Proposition~4.4.4]{kedlaya-part1}, $g^* \calE$ has unipotent local monodromy
along $H$. Since this is true for each $H$,
we may apply the other direction of Theorem~\ref{T:log extend}
 to deduce that $g^* \calE$
is log-extendable on $\overline{Y}$. Hence
$\calE$ admits semistable reduction, as desired.
\end{proof}

Since we are now using a local strategy, it is sensible to refer to
unipotent local monodromy in terms of valuations rather than divisors.
\begin{defn}
Let $\calE$ be an overconvergent $F$-isocrystal on a smooth irreducible
$k$-variety $X$, and let $v$ be a divisorial valuation on $k(X)$ over $k$.
We say that \emph{$\calE$ has unipotent
local monodromy along $v$}
if there exists a birational morphism $f: X' \to X$ such that $v$
is centered on a smooth divisor $Z$ of $X'$, and $f^* \calE$ has
unipotent local monodromy along $Z$. By \cite[Proposition~4.4.1]{kedlaya-part1},
the same will then be true for any other choice of $f$.
Similarly, it is well-defined to say that $\calE$ acquires unipotent local
monodromy along $v$ over a finite separable extension of $k(X)$.
\end{defn}

\subsection{Partial compactifications}

In some applications, it may be helpful to have some sort of semistable
reduction even for isocrystals which are only partially overconvergent.
Here is the correct formulation of the global and local problems.

\begin{defn} \label{D:partial semi}
Let $X$ be a smooth irreducible $k$-variety, and let $\overline{X}$ be
a partial compactification of $X$ (i.e., $\overline{X}$ is a $k$-variety
equipped with an open immersion $X \hookrightarrow \overline{X}$).
Let $\calE$ be an $F$-isocrystal on $X$ overconvergent along
$\overline{X} \setminus X$. We say that $\calE$ \emph{admits semistable
reduction} if there exists a quasiresolution
$(f: X_1 \to \overline{X}, j: X_1 \to \overline{X_1})$ of the pair
$(\overline{X}, \overline{X} \setminus X)$ such that
$f^* \calE$ is log-extendable to $X_1$.
As in Proposition~\ref{P:only perfect}, it is sufficient to check
semistable reduction after extending scalars from $k$ to its algebraic
closure.
\end{defn}

\begin{conj} \label{conj:Shiho refined}
Let $X$ be a smooth irreducible $k$-variety, let $\overline{X}$ be
a partial compactification of $X$, and
let $\calE$ be an $F$-isocrystal on $X$ overconvergent along
$\overline{X} \setminus X$. Then $\calE$ admits semistable reduction.
\end{conj}

\begin{defn} \label{D:local partial semi}
Let $X$ be a smooth irreducible $k$-variety, let $\overline{X}$ be
a partial compactification of $X$, and let $v$ be a valuation on
$k(X)$ centered on $\overline{X}$.
Let $\calE$ be an $F$-isocrystal on $X$ overconvergent along
$\overline{X} \setminus X$. 
We say that $\calE$ \emph{admits local semistable reduction at $v$}
if there exists a quasiresolution
$(f: X_1 \to \overline{X}, j: X_1 \to \overline{X_1})$ of the pair
$(\overline{X}, \overline{X} \setminus X)$ such that
$f^* \calE$ is log-extendable to some open subset of $X_1$ on which
$v$ is centered.
\end{defn}

\begin{remark}
Note that if $\overline{X}$ is proper, then $\calE$ 
is just an overconvergent $F$-isocrystal, and the two
possible interpretations of semistable reduction
(Definitions~\ref{D:semi} and~\ref{D:partial semi}) are consistent;
similarly, the two possible interpretations of local semistable
reduction (Definitions~\ref{D:local semi} and~\ref{D:local partial semi})
are consistent.
More generally, equivalent partial compactifications in the sense of
\cite[Definition~4.1.2]{kedlaya-part1} give rise to equivalent categories
of $F$-isocrystals, and to equivalent notions of global and local
semistable reduction.
\end{remark}

\begin{prop} \label{P:local to global2}
Let $X$ be a smooth irreducible $k$-variety, and let $\overline{X}$ be
a partial compactification of $X$.
Let $\calE$ be an $F$-isocrystal on $X$ overconvergent along
$\overline{X} \setminus X$. 
Suppose that $\calE$ admits local semistable reduction at every valuation 
on $k(X)$ over $k$ centered on $\overline{X}$.
Then $\calE$ admits semistable reduction.
\end{prop}
\begin{proof}
As in Proposition~\ref{P:local to global}.
\end{proof}

\section{Simplification of the local problem}
\label{sec:reduction}

In this chapter, we demonstrate that the local semistable reduction problem
need only be considered around 
valuations of height 1, by an inductive argument.
We also show that valuations whose residue fields have positive transcendence
degree over the base field need not be treated separately, by comparison
between an isocrystal and its ``generic fibre''.

\subsection{\'Etale covers of affine spaces}

Besides de Jong's alterations theorem, it will also be useful to have a 
method for
pushing forward isocrystals onto simple spaces. The following result 
\cite[Theorem~2]{kedlaya-etale2} (based on a technique of Abhyankar for
constructing
finite \'etale morphisms in positive characteristic)
will be of use in this regard.

\begin{prop} \label{P:finite etale}
Let $X$ be an irreducible $k$-variety of dimension $n$, 
let $x \in X$ be a smooth point (whose
existence forces $X$ to be geometrically reduced), and let $D_1, \dots, D_m$ be smooth
irreducible divisors in $X$ meeting transversely at $x$.
Then there exists an open neighborhood $U$ of $x$ in $X$ and a finite \'etale
morphism $f: U \to \AAA^n_k$ such that $D_1, \dots, D_m$ map to coordinate
hyperplanes.
\end{prop}

The relevance of Proposition~\ref{P:finite etale} to our study comes
from the following observation.
\begin{lemma} \label{L:pushforward}
Let $(X,Z)$ be a smooth pair over $k$, put $U = X \setminus Z$, let
$f: Y \to X$ be a finite \'etale morphism of $k$-varieties, 
and let $\calE$ be an isocrystal on $f^{-1}(U)$ overconvergent
along $f^{-1}(Z)$. Then
$\calE$ is log-extendable to $Y$ if and only if $f_* \calE$ is 
log-extendable to $X$.
\end{lemma}
\begin{proof}
By \cite[Theorem~6.4.5]{kedlaya-part1}, we may check log-extendability
by checking unipotence along each component of the boundary divisor.
On one hand, if $f_* \calE$ has unipotent local monodromy,
then so is $f^* f_* \calE$;
however, $\calE$ injects into $f^* f_* \calE$ by adjunction
(see \cite[Definition~2.6.8]{kedlaya-part1}), so by
\cite[Proposition~3.2.20]{kedlaya-part1} it too has unipotent local monodromy. 
On the other hand, suppose $\calE$ has unipotent local monodromy; 
we can then push forward
a log-extension of $\calE$ to obtain a log-extension of $f_* \calE$.
\end{proof}

\subsection{Composite valuations}

We next show that the semistable reduction problem can
be reduced to the restricted local semistable reduction problem which is
only centered at valuations of height 1. 

We first formulate an extension of the full faithfulness theorem for
overconvergent-to-convergent restriction \cite[Theorem~1.1]{kedlaya-full}.
Although we only need the case of a smooth pair, for future
reference we formulate
the general theorem and prove it using a descent argument.
\begin{theorem} \label{T:full}
Let $U \hookrightarrow X$ be an open immersion of
$k$-varieties with dense image, with $U$ smooth.
Then the restriction functor from $F$-$\Isoc^\dagger(U,X/K)$ to
$F$-$\Isoc(U/K)$ is fully faithful.
\end{theorem}
\begin{proof}
Let $\Hom(\calO_U,\calE; U,X/K)$ and
$\Hom(\calO_U,\calE; U/K)$ be the morphisms from $\calO_U$ to $\calE$
in the categories $F$-$\Isoc^\dagger(U,X/K)$ and
$F$-$\Isoc(U/K)$. (These morphisms can be identified with
$F$-invariant horizontal sections on appropriate realizations of $\calE$.)
For any $\calE \in
F$-$\Isoc^\dagger(U,X/K)$, restriction induces an injection
$\Hom(\calO_U,\calE; U,X/K) \to \Hom(\calO_U,\calE; U/K)$,
and the desired result is that this arrow is always surjective.

First suppose $(X, X \setminus U)$ is a smooth pair (this is the only case
that will be used in this paper).
Then this statement follows from \cite[Proposition~6.2.1]{tsuzuki-duke}, under
the assumption of \cite[Conjecture~2.3.2]{tsuzuki-duke}.
However, the latter conjecture is verified by
\cite[Theorem~5.1]{kedlaya-full}, so we may unconditionally deduce the desired
result.

In the general case, we may assume $X$ is irreducible.
Choose a quasiresolution
$(f: X_1 \to X, j: X_1 \hookrightarrow \overline{X_1})$
of the pair $(X, X \setminus U)$,
and put $U_1 = f^{-1}(U)$. Then $(U_1, X_1)$ is a smooth pair;
hence given
$\calE \in F$-$\Isoc^\dagger(U,X/K)$, we may apply
the previous paragraph to show that 
the map 
$\Hom(\calO_{U_1},f^* \calE; U_1,X_1/K) \to \Hom(\calO_{U_1},f^*\calE; U_1/K)$
is bijective.

Suppose we are given $\bv \in \Hom(\calO_U,\calE; U/K)$; we can
pull $\bv$ back to $\Hom(\calO_{U_1},f^*\calE; U_1/K)$. By the previous
paragraph, this element lifts to $\Hom(\calO_{U_1},f^* \calE; U_1,X_1/K)$.
Let $U'$ be the subscheme of $U$ over which $f$ is finite \'etale,
and put $U'_1 = f^{-1}(U')$. We can restrict $\bv$ to
$\Hom(\calO_{U'_1},f^* \calE; U'_1,X_1/K)$, which by
adjunction for finite \'etale morphisms
\cite[Definition~2.6.8]{kedlaya-part1},
\cite[\S 5.1]{tsuzuki-duke} is equal to
$\Hom(\calO_{U'},\calE; U',X/K)$. (Note that we are overloading notation
slightly, by using $\calE$ and $f^* \calE$ to refer also to the restrictions
to $F$-$\Isoc^\dagger(U', X/K)$ and $F$-$\Isoc^\dagger(U'_1,X_1/K)$,
respectively.)

By \cite[Theorem~5.2.1]{kedlaya-part1},
$\bv$ lifts to a morphism in $\Isoc^\dagger(U,X/K)$.
Since we can check compatibility with Frobenius over $U'$,
we have $\bv \in F$-$\Isoc^\dagger(U,X/K)$, as desired.
(A proof using Tsuzuki's cohomological
descent theorem \cite[Theorem~2.1.3]{tsuzuki-descent} is also possible.)
\end{proof}

We next verify a particular geometric instance of the general statement
we are after; ultimately we will reduce back to this case.

\begin{lemma} \label{L:local composite}
Put $\AAA^m_k = \Spec k[t_1, \dots, t_m]$ and put $D = \Spec V(t_1\cdots t_m)
\subset \AAA^m_k$.
Put $\AAA^n_k = \Spec k[u_1, \dots, u_n]$ and put $E = \Spec V(u_1\cdots u_n)
\subset \AAA^n_k$.
Let $\calE$ be an $F$-isocrystal on $(\AAA^m_k \setminus D)
\times (\AAA^n_k \setminus E)$ 
overconvergent along $D \times (\AAA^n_k \setminus E) \cup (\AAA^m_k \setminus D)
\times E$,
with unipotent monodromy along each component
of $D \times (\AAA^n_k \setminus E)$.
Apply Theorem~\ref{T:log extend} to extend $\calE$ to a convergent
log-isocrystal on $(\AAA^m_k \times (\AAA^n_k \setminus E), D \times
(\AAA^n_k \setminus E))$, then restrict to $\{0\} \times (\AAA^n_k \setminus
E)$. Let $\calF$ be the resulting convergent $F$-isocrystal,
and suppose that 
$\calF$ is log-extendable to $\{0\} \times \AAA^n_k$.
Then for any sufficiently large integer $N$, $\calE$ is log-extendable to
\[
\Spec k[t_1 (u_1 \cdots u_n)^{-N}, \dots, t_m (u_1 \cdots u_n)^{-N}, u_1,
\dots, u_n].
\]
\end{lemma}
\begin{proof}
To check log-extendability of $\calE$ to
\[
\Spec k[t_1 (u_1 \cdots u_n)^{-N}, \dots, t_m (u_1 \cdots u_n)^{-N}, u_1,
\dots, u_n],
\]
by Theorem~\ref{T:log extend} it suffices to check unipotence
along $V(u_j)$ for $j=1, \dots, n$.
By generization in the sense of \cite[Proposition~3.4.3]{kedlaya-part1}
(or more precisely, from \cite[Proposition~4.4.1]{kedlaya-part1}),
we may reduce to the case $n=1$.

We may realize $\calE$ as a $\nabla$-module on a space of the form
$A^m_{K,t}[\epsilon, 1] \times A^1_{K,u}[\delta,1]$
for some $\delta, \epsilon \in (0,1) \cap \Gamma^*$. 
By \cite[Lemma~5.1.1(b)]{kedlaya-part1}, for suitable $\delta$,
we can extend $\calE$ to a log-$\nabla$-module $\calE'$
with nilpotent residues on
$A^m_{K,t}[0,1] \times A^1_{K,u}[\delta,1]$.

The restriction of $\calE'$ to $\{0\} \times A^1_{K,u}[1,1]$
is isomorphic to $\calF$, which we assumed admits a log-extension
$\calG$. By Theorem~\ref{T:full}, again for suitable $\delta$
the restriction of $\calE'$ to $\{0\} \times A^1_{K,u}[\delta,1]$ 
is isomorphic to a corresponding restriction of $\calG$.
(Here we are using that the restriction of $\calE'$ is 
overconvergent with respect to $u_1$; this follows from the
same fact on $\calE'$ itself. By \cite[Lemma~3.1.6]{kedlaya-part1}, we may
check this after restriction to the subspace on which
$|t_1| = \cdots = |t_n| = 1$, where
we are given that $\calE$ is overconvergent.)

Pick any $\eta \in (\epsilon,1) \cap \Gamma^*$.
By \cite[Proposition~3.5.3]{kedlaya-part1}, there exists
$\delta \in (0,1) \cap \Gamma^*$ such that $\calE$ is unipotent on
$A^m_{K,t}[\eta,\eta^{1/q}] \times A^1_{K,u}[\delta,1]$.
In other words, on that space, $\calE$ is isomorphic to a successive
extension of $\nabla$-modules pulled back from $A^1_{K,u}[\delta,1]$.
By the hypothesis on the log-extendability of $\calF$,
plus \cite[Proposition~3.6.9]{kedlaya-part1},
we can choose
$\delta$ so that the resulting $\nabla$-modules on $A^1_{K,u}[\delta,1]$
all become unipotent on $A^1_{K,u}[\delta,1)$. Hence $\calE$ admits
a filtration with trivial successive quotients on
\[
\{(t_1, \dots, t_m, u_1) \in \AAA^{m+1}_K: \quad
\delta \leq |u_1| < 1; \quad \eta \leq |t_i| \leq \eta^{1/q}
 \qquad (i=1,\dots,m)\}.
\]
In particular, for any $N$ with $\eta^{1/N} \geq \delta$
(which holds for $N$ sufficiently large), this space contains
\[
\{
(t_1, \dots, t_m, u_1) \in \AAA^{m+1}_K: \quad
\eta^{1/N} \leq |u_1| \leq \eta^{1/(qN)}; \quad 
|t_i/u_1^N| = 1 \qquad (i=1,\dots,m)\},
\]
so $\calE$ is unipotent on the latter. By applying Frobenius
repeatedly, we see that
for each nonnegative integer $h$, $\calE$ is unipotent on
\[
\{
(t_1, \dots, t_m, u_1) \in \AAA^{m+1}_K: \quad
\eta^{1/(q^h N)} \leq |u_1| \leq \eta^{1/(q^{h+1}N)}; \quad 
|t_i/u_1^N| = 1 \qquad (i=1,\dots,m)\}.
\]
By glueing, $\calE$ is unipotent on
\[
\{
(t_1, \dots, t_m, u_1) \in \AAA^{m+1}_K: \quad
\eta^{1/N} \leq |u_1| < 1; \quad
|t_i/u_1^N| = 1 \qquad (i=1,\dots,m)\}.
\]
Hence $\calE$ has unipotent monodromy along the subspace $V(u_1)$ in
$\Spec k[t_1/u_1^N, \dots, t_m/u_1^N, u_1]$,
so Theorem~\ref{T:log extend} yields the desired result.
\end{proof}

We now state a partially restricted version of the
local semistable reduction problem; we will restrict even further
in Conjecture~\ref{conj:local semi}.
\begin{conj} \label{conj:local semi height 1}
Let $X$ be a smooth irreducible $k$-variety, let $\overline{X}$ be
a partial compactification of $X$,
and let $\calE$ be an $F$-isocrystal on $X$ overconvergent along
$\overline{X} \setminus X$.
Then $\calE$ admits local semistable reduction at any
valuation $v$ on $k(X)$ over $k$ of height $1$ centered on $\overline{X}$.
\end{conj}

\begin{prop} \label{P:height 1 to global}
Suppose that Conjecture~\ref{conj:Shiho refined} holds for all varieties of dimension
$<n$, and that Conjecture~\ref{conj:local semi height 1} holds for all 
varieties
of dimension $n$. Then Conjecture~\ref{conj:Shiho refined} also holds for varieties of
dimension $n$. In particular, Conjecture~\ref{conj:local semi height 1}
(for a given $k$ and $K$) implies Conjecture~\ref{conj:Shiho refined} (for the
same $k$ and $K$).
\end{prop}
\begin{proof}[Proof]
Let $X$ be a smooth irreducible $k$-variety of dimension $n$, 
let $\overline{X}$ be
a partial compactification of $X$,
and let $\calE$ be an $F$-isocrystal on $X$ overconvergent along
$\overline{X} \setminus X$.
By Proposition~\ref{P:local to global2}, it suffices to show that for any
valuation $v$ on $k(X)$ of height greater than 1 centered on $\overline{X}$, 
$\calE$ admits local semistable reduction at $v$.
As in Definition~\ref{D:composition}, write $v$ as a composition
$v' \circ \overline{v}$, where $v'$ is a valuation on $X$ of height 1
centered on $\overline{X}$.

We establish a series of reductions of this statement to more restrictive
versions. To begin with, we may assume by
Proposition~\ref{P:only perfect} that:
\begin{enumerate}
\item[(a)] The field $k$ is algebraically closed.
\end{enumerate}
Note that at any point, we may pull back along an alteration and replace $v$
by each of its extensions in turn; Lemma~\ref{L:rel center} guarantees that
these stay centered in the right places. By Theorem~\ref{T:alterations}, we may
thus assume that:
\begin{enumerate}
\item[(b)] There exists a smooth pair $(\overline{Y}, \overline{D})$,
such that $X = Y \setminus D$ for $Y$ the complement of a union of components
of $\overline{D}$ and $D = Y \cap \overline{D}$,
such that $v'$ is centered on $Y$ and $v$ is centered on $\overline{Y}$.
\end{enumerate}
Note also that
the condition that $\calE$ is log-extendable to an open subset of $Y$
on which $v$ is centered is local on $X$ and $Y$, thanks to the 
full faithfulness aspect of Theorem~\ref{T:log extend}.

By the hypothesis that Conjecture~\ref{conj:local semi height 1} holds for all
varieties of dimension $n$, we know that $\calE$ admits local semistable
reduction at $v'$. 
Hence
by passing up a suitable quasiresolution and shrinking,
we may thus assume that:
\begin{enumerate}
\item[(c)] $\calE$ is log-extendable to $Y$.
\end{enumerate}
By shrinking $Y$ and enlarging $D$, we can ensure that:
\begin{enumerate}
\item[(d)] The intersection $E$ of all
of the components of $D$ is nonempty and irreducible, 
and the center of $v'$ on $Y$
is equal to $E$.
\end{enumerate}
By shrinking $X$ and $\overline{Y}$, then
applying Proposition~\ref{P:finite etale} and Lemma~\ref{L:pushforward},
we can ensure that 
\begin{enumerate}
\item[(e)]
We have $Y = \AAA^m_k \times X'$ and
$\overline{Y} = \AAA^m_k \times Y'$, and writing 
$\AAA^m_k = \Spec k[t_1, \dots, t_m]$, we have $D = V(t_1\cdots t_m)$.
\end{enumerate}
By Theorem~\ref{T:alterations} again (applied this time to $(Y', Y' \setminus X')$), 
we may assume that:
\begin{enumerate}
\item[(f)]
There exists a smooth pair $(Y', D')$ with 
$X' = Y' \setminus D'$, such that $\overline{v}$
is centered on $Y'$.
\end{enumerate}
By applying Proposition~\ref{P:finite etale} and Lemma~\ref{L:pushforward},
we can ensure that:
\begin{enumerate}
\item[(g)] We have $Y' = \AAA^{n-m}$, and writing
$\AAA^{n-m}_k = \Spec k[u_1, \dots, u_{n-m}]$, we have
$D' = V(u_1 \cdots u_{n-m})$.
\end{enumerate}
By \cite[Lemma~5.1.1]{kedlaya-part1} (applied on affine subspaces
of $\PP^m_k \times \PP^n_k$),
we may realize the log-extension of $\calE$ to $Y$
as a log-$\nabla$-module with nilpotent residues on $A^m_{K,t}[0, 1] 
\times A^{n-m}_{K,u}[\delta,1]$ for some $\delta 
\in (0,1) \cap \Gamma^*$, which is convergent
with respect to the parameters $t_1, \dots, t_m, u_1, \dots, u_{n-m}$.
The restriction of this log-$\nabla$-module to
$\{0\} \times A^{n-m}_{K,u}[\delta,1]$ represents an 
$F$-isocrystal on $X'$ overconvergent along $D'$.
Let $\calF$ denote the underlying convergent $F$-isocrystal on $X'$;
since 
$\overline{v}$ is a well-defined valuation on $k(X')$
and $\dim(X') < n$, we may invoke the induction hypothesis
in order to ensure that:
\begin{enumerate}
\item[(h)] $\calF$ is log-extendable to $Y'$.
\end{enumerate}
This disturbs restriction (g), but we may
apply Proposition~\ref{P:finite etale} and Lemma~\ref{L:pushforward}
to reestablish it without losing any of the other restrictions.
The desired result in this case follows from Lemma~\ref{L:local composite}.
\end{proof}

\subsection{Positive transcendence degree}

We now give an argument to eliminate the need for separately
treating valuations whose residue fields are not algebraic over $k$.
This again amounts to generization; 
as in the previous section, we calculate in a simple geometric setting
and then reduce the general case back to the simple one.

\begin{lemma} \label{L:generize}
Let $v$ be a valuation on $k(t_1, \dots, t_m, u_1, \dots, u_n)$ over $k$,
with center on $\AAA^{m+n}_k = \AAA^m_k \times \AAA^n_k$
equal to $\{0\} \times \AAA^n_k$; note that this implies that
$v$ is trivial on $\ell = k(u_1, \dots, u_n)$.
Put $D = V(t_1 \cdots t_m) \subset \AAA^m_k$ and $D_\ell = 
V(t_1 \dots t_m) \subset \AAA^m_\ell$.
Let $\calE$ be an $F$-isocrystal on
$(\AAA^m_k \setminus D) \times \AAA^n_k$ overconvergent along
$D \times \AAA^n_k$.
Let $L$ be
the $p$-adic completion of $K(u_1, \dots, u_n)$,
and let $\calF$ be the induced 
isocrystal on $\AAA^m_{\ell} \setminus D_\ell$ overconvergent
along $D_\ell$, with coefficient
field $L$ (as in \cite[Proposition~3.4.3]{kedlaya-part1}).
If $\calF$ admits local semistable reduction at $v$, then
so does $\calE$.
\end{lemma}
\begin{proof}
Choose a quasiresolution $(f_1: X_1 \to (\AAA^m_{\ell} \setminus D_{\ell}),
j_1: X_1 \hookrightarrow \overline{X_1})$ of 
$\AAA^m_{\ell} \setminus D_{\ell}$, such that
$f_1^* \calF$ is log-extendable on an open subset 
on which each extension of $v$ to $\ell(X_1)$ is centered.
Note that each such valuation has center equal to a closed point;
by Theorem~\ref{T:log extend}, $f_1^* \calF$ has unipotent local
monodromy along any divisor passing through that point.

Choose a quasiresolution $(f_2: X_2 \to \AAA^{m+n}_k, 
j_2: X_2 \hookrightarrow \overline{X_2})$ of the pair
$(\AAA^{m+n}_k, D \times \AAA^n_k)$,
such that
$k(X_2)$ contains the maximal separable subextension of
the normal closure of $\ell(X_1)$ over $\ell(t_1, \dots, t_m)
= k(t_1, \dots, t_m, u_1, \dots, u_n)$. 
Put $U = f_2^{-1}((\AAA^m_k \setminus D) \times \AAA^n_k)$,
so that $(\overline{X_2}, \overline{X_2} \setminus U)$ is a smooth
pair and $f_2^* \calE$ is an $F$-isocrystal on $U$
overconvergent along $X_2 \setminus U$.

Let $w$ be any extension of $v$ to $k(X_2)$; we may view $w$
also as an extension of $v$ to $\ell(X_1)$ over $\ell$. 
Let $Y$ be the center of $w$ on $\overline{X_2}$; by
Lemma~\ref{L:rel center}, 
$f_2(Y) = \{0\} \times \AAA^n_k$.
Let $E$ be a component of $\overline{X_2} \setminus U$ 
containing $Y$,
and let $w_E$ be the corresponding
divisorial valuation on $k(X_2)$. Let $v_E$ be the restriction of $w_E$ to
$k(t_1,\dots,t_m,u_1,\dots,u_n)$; then the center of $v_E$ on $\AAA^{m+n}_k$
contains $f_2(Y)= \{0\} \times \AAA^n_k$.

We deduce that $v_E$ is trivial on $\ell$;
we may thus view $w_E$ as a divisorial valuation on $\ell(X_1)$ over $\ell$, whose
center contains the center of $w$. 
As noted above, this means that $f_1^* \calF$ has 
unipotent local monodromy along the center of $w_E$.
By \cite[Proposition~3.4.3]{kedlaya-part1}, $f_2^* \calE$ has unipotent
local monodromy along $E$.

By Theorem~\ref{T:log extend}, we may conclude that $f_2^* \calE$ is 
log-extendable to a subscheme of $\overline{X_2}$ on which each extension
of $v$ to $k(X_2)$ is centered. This implies that $\calE$ admits local
semistable reduction at $v$, as desired.
\end{proof}

\begin{defn}
Let $X$ be an irreducible variety over $k$. By a \emph{minimal valuation}
on $X$, we mean a valuation $v$ on the function field $k(X)$ over $k$
such that $\height(v) = 1$ and $\trdeg(\kappa_v/k) = 0$.
\end{defn}

We now give our most refined version of the local semistable reduction
problem.
\begin{conj} \label{conj:local semi}
Let $X$ be a smooth irreducible $k$-variety, let $\overline{X}$ be
a partial compactification of $X$,
and let $\calE$ be an $F$-isocrystal on $X$ overconvergent along
$\overline{X} \setminus X$. 
Then $\calE$ admits local semistable reduction at any minimal valuation $v$ 
on $k(X)$ centered on $\overline{X}$.
\end{conj}

\begin{theorem} \label{T:minimal to global}
Suppose that for some integer $n$,
Conjecture~\ref{conj:local semi} holds for varieties of dimension
at most $n$ for all algebraically closed $k$.
Then Conjecture~\ref{conj:Shiho refined} holds for varieties of dimension at
most $n$ for all $k$.
\end{theorem}
\begin{proof}
We proceed by induction on $n$; we may thus assume
Conjecture~\ref{conj:Shiho refined} for all varieties of dimension less than $n$.
Let $X$ be a smooth irreducible $k$-variety of dimension $n$, 
let $\overline{X}$ be
a partial compactification of $X$,
and let $\calE$ be an $F$-isocrystal on $X$ overconvergent along
$\overline{X} \setminus X$.
By Proposition~\ref{P:height 1 to global}, it suffices to show that for any
valuation $v$ on $k(X)$ of height 1 centered on $\overline{X}$,
$\calE$ admits local semistable reduction at $v$.
This follows from the assumption of Conjecture~\ref{conj:local semi}
in case $\trdeg(\kappa_v/k) = 0$, so hereafter
we assume instead that $\trdeg(\kappa_v/k) = d > 0$.

As in the proof of Proposition~\ref{P:height 1 to global}, we make
a sequence of reductions, again starting by applying
Proposition~\ref{P:only perfect} to reduce to the case where:
\begin{enumerate}
\item[(a)] The field $k$ is algebraically closed.
\end{enumerate}
By Theorem~\ref{T:alterations}, we may assume that:
\begin{enumerate}
\item[(b)] There exists a smooth pair $(Y,D)$ with
$X = Y \setminus D$, such that $v$ is centered on $Y$.
\end{enumerate}
By Proposition~\ref{P:center dim} (plus Theorem~\ref{T:alterations}
again), we can blow up $X$ and $Y$ to ensure that:
\begin{enumerate}
\item[(c)] The dimension of the center of $v$ on $Y$ is equal to $d$.
\end{enumerate}
By shrinking $X$ and $Y$, we may assume that:
\begin{enumerate}
\item[(d)] $D$ consists of $n-d$ components whose intersection $E$ is the
center of $v$ on $Y$.
\end{enumerate}
By Proposition~\ref{P:finite etale} and Lemma~\ref{L:pushforward},
we may assume that:
\begin{enumerate}
\item[(e)] $Y = \AAA^n_k = \Spec k[t_1, \dots, t_n]$ and
$D = V(t_1\cdots t_{n-d})$.
\end{enumerate}
The desired result now follows from Lemma~\ref{L:generize}.
\end{proof}

\begin{remark}
Theorem~\ref{T:minimal to global} and the $p$-adic local monodromy
theorem imply that local semistable reduction holds at any divisorial
valuation. One way to interpret Theorem~\ref{T:log extend} is that
local semistable reduction at a general valuation $v$ is equivalent to
\emph{uniform} local semistable reduction at all divisorial valuations in
some neighborhood of $v$.
\end{remark}

\begin{remark}
One can deduce refinements of Theorem~\ref{T:minimal to global}
by inspecting its proof and the proof of Proposition~\ref{P:height 1 to global}.
For instance, local semistable reduction for all Abhyankar valuations
follows from local semistable reduction for all monomial valuations.
\end{remark}


\begin{thebibliography}{99}
\bibitem{andre}
Y. Andr\'e, Filtrations de type Hasse-Arf et monodromie $p$-adique,
\textit{Invent. Math.} \textbf{148} (2002), 285--317.

\bibitem{dejong}
A.J. de Jong, Smoothness, semi-stability and alterations,
\textit{Inst. Hautes \'Etudes Sci. Publ. Math.}
\textbf{83} (1996), 51--93.

\bibitem{kedlaya-semicurve}
K.S. Kedlaya, Semistable reduction for overconvergent $F$-isocrystals on a
            curve,
\textit{Math. Res. Lett.} \textbf{10} (2003), 151--159.

\bibitem{kedlaya-local}
K.S. Kedlaya, A $p$-adic local monodromy theorem,
\textit{Annals of Math.} \textbf{160} (2004), 93--184.

\bibitem{kedlaya-full}
K.S. Kedlaya, Full faithfulness for overconvergent $F$-isocrystals,
in \textit{Geometric aspects of Dwork theory},
de Gruyter, Berlin, 2004, 819--835.

\bibitem{kedlaya-etale2}
K.S. Kedlaya,
More \'etale covers of affine spaces in positive characteristic,
\textit{J. Alg. Geom.} \textbf{14} (2005),
187--192.

\bibitem{kedlaya-part1}
K.S. Kedlaya, Semistable reduction for overconvergent $F$-isocrystals, I:
    Unipotence and logarithmic extensions,
\textit{Compos. Math.}, to appear;
arXiv preprint \texttt{math/0405069v5} (2007).

\bibitem{mebkhout}
Z. Mebkhout, Analogue $p$-adique du th\'eor\`eme de Turrittin et le
            th\'eor\`eme de la monodromie $p$-adique,
\textit{Invent. Math.} \textbf{148} (2002), 319--351.

\bibitem{shiho1}
A. Shiho, Crystalline fundamental groups. I. Isocrystals on log
            crystalline site and log convergent site,
\textit{J. Math. Sci. Univ. Tokyo} \textbf{7} (2000), 509--656.

\bibitem{shiho2}
A. Shiho, Crystalline fundamental groups. II. Log convergent cohomology
            and rigid cohomology,
\textit{J. Math. Sci. Univ. Tokyo} \textbf{9} (2002), 1--163.

\bibitem{tsuzuki-duke}
N. Tsuzuki, Morphisms of $F$-isocrystals and the finite monodromy theorem
            for unit-root $F$-isocrystals,
\textit{Duke Math. J.} \textbf{111} (2002), 385--418.

\bibitem{tsuzuki-descent}
N. Tsuzuki, Cohomological descent of rigid cohomology for proper
coverings,
\textit{Invent. Math.} \textbf{151} (2003), 101--133.

\bibitem{vaquie}
M. Vaqui\'e, Valuations, in \textit{Resolution of singularities (Obergurgl, 1997)},
Progr. Math. 181, Birkh\"auser, Basel, 2000, 539--590.

\bibitem{zariski-samuel2}
O. Zariski and P. Samuel,
\textit{Commutative algebra. Vol. II},
Graduate Texts in Mathematics, Vol. 29,
Springer-Verlag, New York, 1975.

\end{thebibliography}
\end{document}